\documentclass[10pt,reqno,a4paper]{amsart}
\usepackage[latin1]{inputenc}
\usepackage{amsmath,amssymb,amsthm}
\usepackage{enumerate}
\usepackage[latin1]{inputenc}
\usepackage{xspace}
\theoremstyle{definition}
\newtheorem{definition}{Definition}[section]

\theoremstyle{remark}
\newtheorem{note}[definition]{Remark}

\theoremstyle{plain}
\newtheorem{lemma}[definition]{Lemma}

\newtheorem{theorem}[definition]{Theorem}

\newtheorem{conjecture}[definition]{Conjecture}

\newcommand{\dividesnot}[2]{{#1 \not \vert #2}}
\newcommand{\tdeg}[1]{\vert#1\rvert}
\newcommand{\Nat}{{\mathbb{N}}}
\newcommand{\Z}{{\mathbb{Z}}}

\newcommand{\Q}{{\mathbb{Q}}}
\def\setsuchas#1#2{\left\{\,{#1}\,\vrule\,{#2}\,\right\}}
\newcommand{\set}[1]{{\{#1\}}}

\newcommand{\PolyRing}[4]{{#1 \lbrack {#2}_{#3}, \dots,
    {#2}_{#4} \rbrack} } 
 
\newcommand{\Pxn}[1]{\PolyRing{#1}{x}{1}{n}}

\newcommand{\Kxn}{\Pxn{K}}

\newcommand{\SQ}{\mathfrak{S}}
\newcommand{\ann}{\mathrm{ann}}

\newcommand{\principalHilb}{p_{n,d}(t)}
\newcommand{\quotientHilb}{q_{n,d}(t)}
\newcommand{\annihilatorHilb}{a_{n,d}(t)}

\newcommand{\qHilb}{q_{n}(t)}
\newcommand{\aHilb}{a_{n}(t)}

\newcommand{\row}{\mathcal{R}}
\newcommand{\col}{\mathcal{C}}

\begin{document}

\newcommand{\anncubdiff}{
\begin{table}[htbp]
  \begin{center}
    \label{tab:anndeg3diff}
    \begin{tabular}{|c|l|} \hline
n & \(q_n(t) - w_n(t)\) \\ \hline
3 & \( 3t(1+t) \) \\ 
4 & \( 3t(1+t)^2 \) \\ 
5 & \( t(1+t)(t^2+8t+1) \) \\ 
6 & \( 9t^2(1+t)^2 \) \\ 
7 & \( 27t^3(1+t) \) \\ 
8 & \( 27t^3(1+t)^2 \) \\ 
9 & \( 3t^3(1+t)(t^2+26t+1) \) \\ 
10 & \( 81t^4(1+t)^2 \) \\ 
11 & \( 243t^5(1+t) \) \\ 
12 & \( 243t^5(1+t)^2 \) \\ 
13 & \( t^5(1+t)(t^2+728t+1) \) \\ 
14 & \( 729t^6(1+t)^2 \) \\ 
15 & \( 2187t^7(1+t) \) \\ 
16 & \( 2187t^7(1+t)^2 \) \\ \hline
    \end{tabular}
    \caption{\(\annihilatorHilb - \principalHilb\) for a cubic generic form}
  \end{center}
\end{table}
}

\newcommand{\taborder}{
\begin{table}[htbp]
  \begin{center}
    \small{
\begin{tabular}{|c|cccccccc|}
\hline
  d & 3 & 5 & 7 & 9 & 11 & 13 & 15 & 17 \\ \hline
n \\
3 & 1\\
4 & 1 \\
5 & 1 & 1 \\
6 & 2 & 1 \\
7 & 3 & 1 & 1 \\
8 & 3 & 2 & 1 \\
9 &3 & 3 &  1 & 1 \\
10 & 4& 3 & 2 & 1 \\
11& 5& 3& 3& 1& 1 \\
12& 5& 4& 3& 2& 1 \\
13& 5& 5& 4& 3& 1& 1 \\
14& 6& 5& 4& 3& 2& 1 \\
15& 7& 6& 5& 3& 3& 1& 1 \\
16& 7& 6& 5& 4& 3& 2& 1 \\
17& -& -& -& 5& 4& 3& 1& 1 \\
18& -& -& -& -& 4& 3& 2& 1 \\
19& -& -& -& -& -& 3& 3& 1\\
20& -& -& -& -& -& -& 3& 2\\ \hline
\end{tabular}
    \caption{Order of \(\annihilatorHilb - \principalHilb\) for small \(n,d\)}
    }
    \label{tab:order}
  \end{center}
\end{table}
}

\newcommand{\tabodd}{
    \begin{table}[htbp]
      \begin{center}
        \caption{Difference between true and predicted Hilbert series of the
          annihilator of a generic form of odd degree}
        \label{tab:odd}
\begin{tabular}{|c|l|c|c|c|c|c|c|c|c|}
\hline 
n & deg=3 & 5 & 7 & 9 & 11 & 13 & 15 & 17 & 19 \\ \hline
3 &0 & & & & & & & & \\ 
4 &0 & & & & & & & & \\ 
5 &t & 0 & & & & & & & \\ 
6 &0 & 0 & & & & & & & \\ 
7 &0 & t & 0 & & & & & & \\ 
8 &0 & 0 & 0 & & & & & & \\ 
9 &\(3t^3\) & 0 & \(t\) & 0 & & & & & \\ 
10 &\(t^4\) & 0 & 0 & 0 & & & & & \\ 
11 &\(t^5\) & \(t^3\) & 0 & \(t\) & 0 & & & & \\ 
12 &\(t^6+12t^5\) & 0 & 0 & 0 & 0 & & & & \\ 
13 &\(t^7+13t^6+t^5\) & 0 & 0 & 0 & \(t\) & 0 & & & \\ 
14 &\(t^8+14t^7+91t^6\) & 0 & 0 & 0 & 0 & 0 & & & \\ 
15 &\(15t^8+105t^7\) & 0 & 0 & \(t^3\) & 0 & t & 0 & & \\ 
16 &\(16t^9+120t^8+559t^7\) & \(t^6\) & 0 & 0 & 0 & 0 & 0 & & \\ 
17 && & & 0 & 0 & 0 & \(t\) & 0 & \\ 
18 && & & & 0 & 0 & 0 & 0 & \\ 
19 && & & & & \(t^3\) & 0 & \(t\) & 0 \\ 
20 && & & & & & 0 & 0 & 0 \\ 
21 && & & & & &   & 0 & \(t\)  \\ \hline
\end{tabular}
      \end{center}
    \end{table}
}

\newcommand{\taboddmac}{
    \begin{table}[htbp]
      \begin{center}
        \caption{Difference between true and predicted Hilbert series of the
          annihilator of a generic form of odd degree \(>3\)}
        \label{tab:oddmac}
\begin{tabular}{|c|l|c|c|c|c|c|c|c|}
\hline 
\(n -d\) & deg=5 & 7 & 9 & 11 & 13 & 15 & 17 & 19 \\ \hline
0   & 0     & 0   & 0     & 0   & 0     & 0   & 0   & 0    \\
1   & 0     & 0   & 0     & 0   & 0     & 0   & 0   & 0    \\ 
2   &\(t\)  &\(t\)&\(t\)  &\(t\)& \(t\) &\(t\)&\(t\)&\(t\) \\ 
3   & 0     & 0   & 0     & 0   & 0     & 0   & 0          \\
4   & 0     & 0   & 0     & 0   & 0     & 0   & 0          \\
5   & 0     & 0   & 0     & 0   & 0     & 0   & 0          \\
6   &\(t^3\)& 0   &\(t^3\)& 0   &\(t^3\)                   \\
7   & 0     & 0   & 0     & 0                              \\
8   & 0     & 0   & 0                                      \\
9   & 0     & 0                                            \\
10  & 0                                                    \\
11  &\(t^6\)                                               \\
\hline
\end{tabular}
      \end{center}
    \end{table}
}

  \newcommand{\tabqgens}{
  \begin{table}[htb]
    \begin{center}
      \leavevmode
    \begin{tabular}{|c|c|c|c|c|c|c|c|c|c|c|c|c|} \hline
      n &    2 & 3 & 4 & 5       & 6 & 7       &    8  
      &  9    & 10        & 11          & 12       & 13 \\ \hline
      Diff & 0 & 0 & 0 & \(t^3\) & 0 & \(t^4\) & \(t^4\) &
      \(t^5\) & \(10t^5\) & \(t^6+t^5\) &\(64t^6\) & \(t^7 + 13t^6\)
      \\ \hline
      \end{tabular}
    \end{center}
    \caption{Difference between the true Hilbert series and the
      ``anticipated Hilbert series''  for generic ideals generated by
      two quadratic forms}
    \label{tab:2gens}
  \end{table}
}

\sloppy

\author{Guillermo Moreno-Socías}
\address{Laboratoire GAGE, 
  École Polytechnique \\
 91128 Palaiseau Cedex \\France 
}
\email{moreno@gage.polytechnique.fr}
\author{Jan Snellman}
\address{Department of Mathematics\\
Link\"oping University\\
SE-58183 Link\"oping\\
SWEDEN
}

\email{jasne@mai.liu.se}

\thanks{Snellman was supported grants from
  Svenska institutet and by grant n. 231801F from Centre International 
  des Etudiants et Stagiaires while visiting École Polytechnique, and by
  grants from Svenska Institutet and Kungliga Vetenskapsakademin while
  visiting University of Wales, Bangor} 

\title[Hilbert series of generic ideals in
  the exterior algebra]{Some conjectures about the Hilbert series of
    generic ideals in  the exterior algebra}

%\date{\today} 
\subjclass{15A75,13D40} 
\keywords{exterior algebra, generic ideals, Hilbert series}

\begin{abstract}
  We calculate the Hilbert series of a quotient of the exterior
  algebra by a generic form of even degree, and give conjectures about
  the Hilbert series of other generic quotients. 
\end{abstract}

\maketitle

\begin{section}{Introduction}
  In the symmetric algebra \(\Kxn\), the set of Hilbert series coming
  from homogeneous quotients are classified by Macaulays theorem
  \cite{Macaulay:Enum,Ebud:View,Bigatti95}. There is an infinite
  number of possible series, but if we fix positive integers
  \(d_1,\dots,d_r\),  and restrict our study to quotients by
  homogeneous ideals \(I\) of ``type'' or ``numerical character''
  \((d_1,\dots,d_r)\), ie generated by forms of those prescribed
  degrees, then there are only finitely many Hilbert
  series. Furthermore, in the affine space parametrising these
  homogeneous ideals,  there is a
  Zariski-open subset of ideals with the same 
  Hilbert series, and the Hilbert series obtained on this open set is
  minimal \cite{Froeberg:OnHilb, Froeberg:Inequality}.

  Unfortunately, even though we know the set of \emph{all} Hilbert
  series, we do not know what Hilbert series arise from ideals of
  numerical character \((d_1,\dots,d_r)\). In fact, we do not even
  know the ``generic'' series, but it is conjectured
  \cite{Moreno:Revlex, Froeberg:Inequality} that it is 
  \(\left \langle (1-t)^{-n} \prod_{i=1}^r (1-t^{d_i}) \right
  \rangle\); the brackets mean ``truncate before the first
  non-positive coefficient''. 
 
  In the exterior algebra \(\bigwedge V_n\), we also know the set of
  all Hilbert series of homogeneous quotients, by the so-called
  Kruskal-Katona theorem \cite{Kruskal,Katona,CleLi,Aramova:Gotzman}.
  Here, this set is finite, so one would think that it should be easy to
  find the subset of Hilbert series coming from quotients by ideals having
  a prescribed numerical character. In particular, it should be easy
  to find the generic value. However, very little is known. 

  In this article, we give one new result (the series for a quotient
  by \emph{one} form of \emph{even} degree) and several conjectures,
  supported  by extensive computer calculations. 

  It is worthwhile to point out that the problem of determining the
  Hilbert series of quotients by generic \emph{quadratic} forms is especially 
  interesting, since it determines the Koszulness of
  the quadratic algebras in question. We refer to the recent article
  by Fr\"oberg and L\"ofwall \cite{KosLie}.
\end{section}

\begin{section}{Notation}
  Let \(K\) be a field of characteristic 0. Then \(\Q\) is the prime
  subfield of \(K\).  For any positive integer \(n\), let 
  \(V=V_n\) be an \(n\)-dimensional 
  vector space over \(K\), with a distinguished basis
  \(X_n=\set{x_1,\dots,x_n}\). Let \(\Kxn\) denote the symmetric
  algebra on \(V_n\), and let \(\bigwedge V_n\) denote the exterior
  algebra on \(V_n\). We define \(\SQ(V_n)\), the \emph{square-free} algebra on
  \(V_n\), to be the commutative \(K\)-algebra generated by \(X_n\),
  with the relations \(x_i^2 = 0\); in other words, \(\SQ(V_n) =
  \frac{\Kxn}{(x_1^2,\dots,x_n^2)}\). There is an isomorphism of
  graded vector spaces between \(\bigwedge V_n\) and \(\SQ(V_n)\), but 
  they are not isomorphic as \(K\)-algebras, since the exterior
  algebra is skew-commutative and \(\SQ(V_n)\) is commutative.

  We shall need the following operations for formal power series.
  \begin{definition}
    Let \(f(t) = \sum_{i=0}^\infty a_i t^i \in \Z[[t]]\), \(g(t) =
    \sum_{i=0}^\infty b_i t^i \in \Z[[t]]\). We say that \(f \ge g\)
    if \(a_i \ge b_i\) for all \(i\). We define 
    \begin{align*}
      \max(f(t),g(t)) & = \sum_{i=0}^\infty \max(a_i,b_i) \\
      \left \langle f(t) \right \rangle & = \sum_{i=0}^\ell a_i t^i, \qquad
      \ell = \max(\setsuchas{i}{a_j > 0 \text{ for } j \le i})
       \\
      \left \rangle f(t) \right \langle & = \sum_{i=\ell}^\infty a_i
      t^i, \qquad \ell = \min(\setsuchas{i}{a_j > 0 \text{ for } j \ge i})
      \qquad 
    \end{align*}
    We use the conventions \(\max(\emptyset) = -1 = \min(\Nat)\), 
    \(\min(\emptyset)   = +\infty = \max(\Nat)\).
  \end{definition}

  Let \([X_n]\) denote the free abelian monoid on \(X_n\), and denote
  by \(Y_n\) the subset of square-free monomials. Then \(Y_n\) is a
  \(K\)-basis for both \(\bigwedge V_n\) and \(\SQ(V_n)\). We define
  the \emph{degree}  of a monomial in \([X_n]\) (and in \(Y_n\)) in
  the usual way, and denote by \([X_n]^d\) and \(Y_n^d\) the subset of 
  monomials (square-free monomials) of degree \(d\).

  A form
  \begin{math}
    \bigwedge \Kxn \ni f = \sum_{m \in [X_n]^d} c_m m
  \end{math}
  is said to be \emph{generic} if the coefficients \(c_m \in K\)
  fulfil the following conditions:
  \begin{enumerate}
  \item \(c_m \not \in \Q\),
  \item \(m \neq m' \implies c_m \neq c_{m'}\),
  \item The set of all \(c_m\)'s is algebraically independent over
    \(\Q\). 
  \end{enumerate}
  A homogeneous  ideal \(I \subset \Kxn \) is called generic if it can be
  minimally generated by a finite set of generic forms, so that all of 
  the occuring 
  coefficients of the forms are different, and so that the set of all
  occuring coefficients is algebraically independent over \(\Q\).
  If the forms have degrees \(d_1,\dots,d_r\), then we say that \(I\) has 
``numerical character'' \((d_1,\dots,d_r)\).
  It is an important fact that  any two generic   
  ideals  of 
  the same numerical character
  have  
  the same initial ideal and the same Hilbert series.  

  Now consider the affine space \(V=\mathbf{A}^{\binom{n+d_1 -1}{d_1}}
    \times \cdots \times \mathbf{A}^{\binom{n+d_r -1}{d_r}}\)
    parametrising the set  
  of homogeneous ideals of numerical character \((d_1,\dots,d_r)\). 
  Since there are
  countably many conditions to be fulfilled for an ideal to be
  generic, the subset of the 
  parameter space corresponding to generic ideals is not open, but a
  countable intersection of open sets, hence dense.
  However, in \(V\) there is a 
  Zariski-open subset corresponding to ideals with the same 
  Hilbert function, and the generic ideals are contained in this
  subset \cite{Froeberg:OnHilb}. 

  We make similar definitions for the square-free algebra, and for
  the exterior algebra. Here, a generic form is  a generic linear combination of 
\emph{square-free} monomials of a certain degree. It is still true
that the generic Hilbert series is attained on an open component of
the parameter space, and that the generic ideals are contained in this 
component.

\end{section}

\begin{section}{Hilbert series for generic principal ideals in the
    symmetric and square-free algebra} 
  \begin{subsection}{Principal ideals in the symmetric algebra}
    If \(f \in \Kxn\) is a non-zero form of degree \(d\), not
    necessarily homogeneous, then clearly the Hilbert series of the
    quotient \(\frac{\Kxn}{(f)}\) is \((1-t)^{-n} (1-t^d)\). 
\end{subsection}

\begin{subsection}{Principal ideals in the square-free algebra}
If \(f \in
    \SQ(V_n)\) is a generic form of degree \(d\), then there is a
    similar simple formula for \(\frac{\SQ(V_n)}{(f)}(t)\) (the
    Hilbert series of the quotient). To state the formula, we need
    some additional notation.

  \begin{definition}
        We denote the zero series by \(0\), and define
    \begin{align*}
      \Delta_{n,d}(t) &=     \left \rangle (t^d -1) (1+t)^n \right \langle  \\
    &=\sum_{v=\lceil (n-d)/2 \rceil}^n
    \left(\binom{n}{v} - \binom{n}{v+d}\right) t^v \\
    \delta_{n,d}(t) &= \left \langle (1+t)^n (1-t^{d}) \right \rangle \\
    &=  \sum_{v=0}^{\lfloor (n-d)/2 \rfloor}  (\binom{n}{v+d}-
    \binom{n}{v}) t^v  
    \end{align*}
  \end{definition}

  The following result is due to Fröberg \cite{Froeberg:HilbGenForm}.
  \begin{theorem}\label{thm:principal}
    Let \(f \subset \SQ(V_n)\) be a generic form of degree \(d\).
    Then 
    \begin{equation}
      \label{eq:SQser}
      \frac{\SQ(V_n)}{(f)}(t) = \delta_{n,d}(t)
    \end{equation}
  \end{theorem}
  \begin{proof}
    By considering the graded exact sequence 
    \begin{equation}
      \label{eq:exseqSQ}
      0 \longrightarrow \ann(f)(-d) \longrightarrow \SQ(V_n)(-d)
      \xrightarrow{\cdot f} \SQ(V_n) \longrightarrow
      \frac{\SQ(V_n)}{(f)} \longrightarrow 0   
    \end{equation}
    in each degree \(r\),
    we see that \eqref{eq:SQser} holds if and only if multiplication by
    \(f\), regarded as a linear map \(\phi_r\) from \(\SQ(V_n)_r\) to
    \(\SQ(V_n)_{r+d}\), is injective when \(\binom{n}{r} \le
    \binom{n}{r+d}\), and surjective when \(\binom{n}{r} \ge
    \binom{n}{r+d}\). 

    Write \(f=\sum_{m \in Y_n^d} c_m m\).
    For \(0 \le r \le n-d\),  \(Y_n^r\) is a basis of \(\SQ(V_n)_r\), and
    \(Y_n^{r+d}\) is a basis of \(\SQ(V_n)_{r+d}\). Thus, we must show 
    that for each \(r\), the matrix of \(\phi_r\) in this basis has
    maximal rank. This matrix has rows indexed by \(Y_n^{r+d}\)
     and columns indexed by \(Y_n^{r}\).
    The entry at row \(R\), column \(C\) is 
    \begin{displaymath}
      \begin{cases}
        0 &  \dividesnot{C}{R}\\
        c_m & R=mC
      \end{cases}
    \end{displaymath}
    If we specialise this matrix, the rank can only decrease, so if we 
    can prove that some specialised matrix has full rank, then we are
    done. Putting all \(c_m=1\), we obtain the incidence matrix of
    \(r\)-subsets of \([n]\) into \(r+d\)-subsets of \([n]\), that is,
    the rows are indexed by \(r\)-subsets and the columns by
    \(r+d\)-subsets, with a \(1\) at the \(a,b\)'th position iff \(a
    \subset b\), and 0 otherwise. It has been shown by  combinatorialists
    that this matrix has full rank \cite{Wilson:design,
      Kantor:incidence, Graver:design}. 
  \end{proof}
\end{subsection}

\end{section}

\begin{section}{Principal ideals in the exterior  algebra --- the
    difference between even and odd degree}

Let \(f \in \bigwedge V_n\) be a generic form of degree \(d\).
Denote the Hilbert series of \(\frac{\bigwedge V_n}{f}\) by \(\quotientHilb\),
that of the annihilator of \(f\) by \(\annihilatorHilb\), and that of the
principal ideal \((f)\) by \(\principalHilb\). From the
the graded  exact sequence
\begin{equation}
  \label{eq:exseqWEDGE}
  0 \longrightarrow \ann(f)(-d) \longrightarrow \bigwedge V_n(-d)
  \xrightarrow{\cdot f} \bigwedge V_n  \longrightarrow
  \frac{\bigwedge V_n}{(f)} \longrightarrow 0   
\end{equation}
we get that 
\begin{equation}
  \label{eq:dimsum}
  \begin{split}
    \quotientHilb  & = t^d \annihilatorHilb - t^d(1+t)^n + (1+t)^n \\
    &= t^d \annihilatorHilb + (1+t)^n(1-t^d) \\
    \annihilatorHilb &= t^{-d} \left( \quotientHilb - (1+t)^n(1-t^d) \right)
  \end{split}
\end{equation}

If \(d\) is even, we shall prove that the vector space map 
\begin{equation}
  \label{eq:fmul}
  \bigwedge^v V_n
  \xrightarrow{\cdot f} \bigwedge^{v+d} V_n  
\end{equation}
is injective ``when it can be'', ie when \(\binom{n}{v} \le
\binom{n}{v+d}\), and surjective ``when it can be'', ie when \(\binom{n}{v} \ge
\binom{n}{v+d}\). This leads immediately to the formul{\ae}
\begin{equation}
\begin{split}
\quotientHilb &=\left \langle (1+t)^n (1-t^{d})
\right \rangle =\delta_{n,d}(t) \\
  \annihilatorHilb &= t^{-d} \left( \quotientHilb - (1-t^d)(1+t)^n \right) \\
  &= t^{-d} \left( \delta_{n,d}(t) - (1-t^d)(1+t)^n \right) \\
  & = t^{-d} \sum_{r=0}^n  \Biggl[\max\left(0,\binom{n}{r+d} -
    \binom{n}{r}\right) - \left( 
  \binom{n}{r+d} - \binom{n}{r}\right) \Biggr]  t^r \\
  & = t^{-d} \sum_{r=0}^n \max \left(0, -\binom{n}{r+d} + \binom{n}{r} 
    \right) t^r
  \\
  & = t^{-d} \Delta_{n,d}(t)
    \end{split}
\end{equation}
In particular, as \(n \to \infty\), \((1+t)^{-n}\quotientHilb \to
(1-t^d)\), and \(\annihilatorHilb \to 0\), with respect to the \(t\)-adic
norm on \(\Z[[t]]\).

    If \(d\) is odd, then   we have that \(f^2 = 0\), hence
  \(fg=0\) whenever \(g \in (f)\), hence
  \(\ann(f) \supseteq (f)\), hence \(\annihilatorHilb \ge \principalHilb\).
    In other words,
    there is a graded complex 
    \begin{equation}
      \label{eq:infcomp}
   \left(\bigwedge V\right)(-d) \xrightarrow{\cdot f}
\bigwedge V \xrightarrow{\cdot 
        f} \left(\bigwedge V\right)(d) 
    \end{equation}
    the graded homology of which determines \(\annihilatorHilb -
    \principalHilb\). 
    In the (not very interesting) case \(d=1\), then we know from
    \cite{Aramova:Cohomology} that this homology vanishes.
    For odd \(d > 1\), we  guess
  that for a fixed degree \(r\), and \(n\) very large, this homology
  vanishes. Hence, in degree \(r\), the ``obstruction to injectivity'' 
  in \eqref{eq:fmul}  is as small as possible. An equivalent
  formulation: consider the start of a minimal free graded \(\bigwedge 
  V_n\)-resolution of \(\frac{\bigwedge V_n}{(f)}\),
  \begin{displaymath}
    \frac{\bigwedge V_n}{(f)} \leftarrow \bigwedge V_n
    \xleftarrow{\cdot f} \bigwedge V_n \leftarrow \bigoplus_{j=1}^r
    (\bigwedge V_n)(-\beta_{2,i}),
  \end{displaymath}
  where \(\beta_{2,i}\) are the graded Betti numbers. Then we guess
  that as \(n\) increases, and for a fixed \(i \neq 2d\),
  \(\beta_{2,i} = 0\). On the other hand, for sufficiently large
  \(n\), we guess that \(\beta_{2,2d} = 1\).
  Since \(\beta_{2,i}\) is the dimension of the degree \(i-d\) part of 
  a certain Tor group, this
  conjecture can also be stated in terms of Cartan homology (see
  \cite{Aramova:Gotzman}).

We show the order (ie the smallest \(\ell\) for which \(t^\ell\)
occurs with non-zero coefficient) 
     of  \(\annihilatorHilb -  \principalHilb \) for small \(n,d\) in
%    Table~\ref{tab:order}. 
    Table 1.
    \taborder
It would seem that the order of the
    difference grows linearly in \(n\), so that \(\annihilatorHilb -
    \principalHilb \to 0\) rather rapidly.

  Let us turn to the consequences of this conjecture. We get that
  \(\annihilatorHilb \sim \principalHilb\) with respect to the \((t)\)-adic
  filtration. It then follows from \eqref{eq:dimsum} that 
  \begin{equation}
    \label{eq:oddconj}
    \quotientHilb \sim t^d \principalHilb  + (1+t)^n(1-t^d) 
  \end{equation}
    Substituting \(\principalHilb = (1+t)^n - \quotientHilb\) and solving for
    \(\quotientHilb\) we get that
    \begin{equation}
      \quotientHilb \sim \frac{(1+t)^n t^d + 
        (1+t)^n(1-t^d)}{(1+t^d)} = \frac{(1+t)^n}{(1+t^d)},
    \end{equation}
      hence
    \begin{equation}\label{eq:oddagain}
      \frac{\frac{\bigwedge V_n}{(f)}(t)}{\bigwedge(V_n)(t)}
      =  \frac{\quotientHilb}{(1+t)^n} \to  \frac{1}{1+t^d} 
\qquad \text{ as } n   
      \to \infty.
    \end{equation}

\end{section}

  \begin{section}{Principal ideals on generic forms of even degree in
      the exterior algebra}
  If \(d=2\) then
  we can change coordinates on \(V\) and replace \(f\)  with 
  the form \(x_1x_2 + x_3x_4 + \cdots\), as is demonstrated in
  \cite{Bourbaki:FormQ}. The Hilbert series of the quotient can now be
  easily calculated. We get that 
  \(\frac{\bigwedge(V_n)}{(f)}(t) = \left \langle (1+t)^n (1-t^2)
  \right \rangle\), which is the same as the Hilbert series for the
  corresponding quotient in the square-free algebra. 

   \begin{note}
     It is \emph{not true} that if \(f_e = \sum_{1\le i<j \le n}
     \alpha_{ij} x_i x_j\) is a non-generic quadratic form in
     \(\bigwedge V_n\),  
     and \(f_s = \sum_{1\le i<j \le n} \alpha_{ij} x_i x_j\) is
     the corresponding form in \(\SQ(V_n)\), then \(\frac{\bigwedge
       V_n}{(f_e)}\) and \(\frac{\SQ(V_n)}{(f_s)}\) have the same
     Hilbert series. For an example, consider the form
     \(x_1x_2 + x_1x_3 + x_1x_4 + x_3x_4\). The quotient of
     \(\bigwedge V_4\) by this form has Hilbert series 
     \begin{math}
       5t^2 + 4t+1,
     \end{math}
     but the corresponding quotient of \(\SQ(V_4)\) has series
     \begin{math}
       t^3+5t^2 + 4t+1. 
     \end{math}
   \end{note}

  We next show that if the degree \(d\) of \(f\) is even, then
  the Hilbert series of the quotient \(\frac{\bigwedge V_n}{(f)}\) is
  the same as for the 
  square-free algebra. To this end, we need some
  combinatorial results, which we have collected in the appendix.
  With the aid of these, we can prove:

\begin{theorem}\label{thm:evenform}
  Let \(f \in \wedge^d V\), with \(d\) even, be a generic form. Then
  the linear transformation  
  \begin{equation}
    \label{eq:multbyf}
    \wedge^r V \xrightarrow{f \cdot} \wedge^{r+d} V
  \end{equation}
  is injective for \(2r+d \le n\), and surjective for \(2r+d \ge n\).
\end{theorem}
\begin{proof}
  We put \(k=r+d\).
  Suppose that
  \begin{equation}
    \label{eq:fis}
    f= \sum_{K \in \binom{[n]}{d}} c_K x_K
  \end{equation}
  The matrix of the map
  \eqref{eq:multbyf} is an \(\binom{n}{r+d} \times \binom{n}{d}\)
  matrix, \(\widetilde{M_{r,r+d,n}}\), where the rows are indexed by
  \((r+d)\)-subsets \(K\), and 
  the columns by \(d\)-subsets \(T\). The entry at position \((K,T)\)
  is 
  \begin{equation}
    \label{eq:pos}
    \begin{cases}
      0 & \text{ if } T \not \subseteq K \\
      \sigma(T, K) c_T& \text{ if } T \subseteq K
    \end{cases}
  \end{equation}
  We must prove that this matrix has
  maximal rank. Clearly, the rank can not increase under
  specialisation, so if we prove that the matrix obtained by replacing 
  each \(c_T\) with 1 has maximal rank, then so does
  \(\widetilde{M_{r,r+d,n}}\). 
  However, the specialised matrix is nothing but the matrix
  \(M_{r,r+d,n}\) 
  of Theorem~\ref{thm:fullrank}, so it has full rank.
  \end{proof}

  \begin{theorem}\label{thm:evenser}
    Let \(f \in \bigwedge V_n\) be a generic form of degree \(d\),
    with \(d\) even. Then 
    \begin{equation}
      \label{eq:evendegHser}
      \frac{\bigwedge V_n}{(f)}(t) = \left \langle (1+t)^n (1-t^{d})
      \right \rangle = \delta_{n,d}(t)
    \end{equation}
  \end{theorem}
  \begin{proof}
    This follows from Theorem~\ref{thm:evenform}, together with
    \eqref{eq:exseqWEDGE}.  
  \end{proof}

\end{section}

\begin{section}{Principal ideals on generic forms of odd degree in
      the exterior algebra}
Let \(d\) be an odd integer.
Recall that we've conjectured that 
    \(\annihilatorHilb -  \principalHilb \to 0\) as \(n \to \infty\),
    and that this 
    conjecture leads to the conclusions that \(p_{n,d}(t) \sim
    (1+t)^n(1+t^d)^{-1}\). In this section, we shall try to guess
    the exact value of \(\quotientHilb\).

Since \(\annihilatorHilb \ge \principalHilb\), \(\annihilatorHilb \ge
\Delta_{n,d}(t)\), it follows that \(\annihilatorHilb \ge
\max(\principalHilb,\Delta_{n,d}(t))\). We tabulate the difference \(\annihilatorHilb -
\max(\principalHilb,\Delta_{n,d}(t))\) in Table~\ref{tab:odd} and
Table~\ref{tab:oddmac}. 
\tabodd

Using the data of Table~\ref{tab:oddmac}, we make the following
conjecture:
    \begin{conjecture}\label{conj:oddfive}
      Let \(d\) be an odd integer \(> 3\). Then, putting
      \(\tau_{n,d}(t)= \annihilatorHilb-\max \left( \principalHilb,
        \Delta_{n,d}(t) 
      \right)\),  
      \begin{equation}
    \tau_{n,d}(t)=
    \begin{cases}
      t^{v(v-1)/2} & \exists v,s \in \Nat: \, v > 0, \, n-d = -1 +
      \frac{5}{2}v  + \frac{1}{2}v^2, \, d = 5 + 2vs \\
      0 & \text{ otherwise}
    \end{cases}
      \end{equation}
    \end{conjecture}

\taboddmac

    This conjecture yields a formula for the Hilbert series, but since 
    said formula is very complicated, we do not write it down; instead
    we show how to derive \(\quotientHilb\). From
    \begin{equation}
      \label{eq:solveme}
      \begin{split}
        \annihilatorHilb &= \tau_{n,d}(t) + \max \left(
          \principalHilb, \Delta_{n,d}(t) 
        \right) \\ 
        \quotientHilb & = \annihilatorHilb t^d + (1+t)^n(1-t^d) \\
        \principalHilb & = (1+t)^n - \quotientHilb
      \end{split}
    \end{equation}
    we get
    \begin{equation}
      \label{eq:sol}
      \begin{split}
        \principalHilb 
        &= (1+t)^n - \quotientHilb \\
        &= (1+t)^n - \annihilatorHilb t^d - (1+t)^n(1-t^d) \\
        &= (1+t)^n - t^d\tau_{n,d}(t) - t^d\max \left( \principalHilb,
          \Delta_{n,d}(t) \right) 
           - (1+t)^n(1-t^d) \\
        &= t^d \left( (1+t)^n - \tau_{n,d}(t)- \max(\principalHilb,
          \Delta_{n,d}(t)) \right)         
      \end{split}
    \end{equation}
    Hence, writing
    \(\principalHilb = \sum_{i=0}^n a_i t^i\), with the \(a_i\)'s as
    undetermined coefficients, and denoting the
    \(t^i\)-coefficient of \(\tau_{n,d}(t)\) by \(b_i\), we get the equation
    \begin{equation}
      \label{eq:iter}
      a_\ell = \binom{n}{\ell - d} - b_{i-\ell} -\max(a_{\ell - d},
      \binom{n}{\ell - d} - 
      \binom{n}{\ell}) 
    \end{equation}
    which we can solve recursively, using the initial values 
    \begin{displaymath}
      a_0 = \cdots = a_{d-1} = 0, \qquad a_d = a_n = 1.
    \end{displaymath}

    For the case \(d=3\), we proceed differently: we tabulate
    \(q_{n,3}(t) - w_{n,3}(t)\) in Table 4,
    %~\ref{tab:anndeg3diff},
    and from that, make the following conjecture:
    \begin{conjecture}\label{conj:deg3}
      The Hilbert series of \(\frac{\bigwedge V_n}{(f)}\), where \(f\) is a
      generic cubic form, is given by
      \begin{equation}
        \label{eq:pdeg3}
        \begin{split}
          p_{n,3}(t) &= \frac{t^d L_n(t) + (1+t)^n}{1+t^d} \\         
      L_n(t) &=
      \begin{cases}
        (3t)^{2\ell -1} (1+t)^2 & n = 4 \ell \\
        c_1(n)t^{2\ell -1}(1+t)(1+(3^{c_2(n)}-1)t + t^2) & n = 4 \ell + 1\\
        (3t)^{2\ell} (1+t)^2 & n = 4 \ell + 2\\
        (3t)^{2\ell +1} (1+t) & n = 4 \ell + 3
      \end{cases}
        \end{split}
      \end{equation}
      where \(c_1(n),c_2(n)\) are some positive integers.
    \end{conjecture}
    \anncubdiff

\end{section}

\begin{section}{Hilbert series for generic non-principal ideals in the
    symmetric and square-free algebra} 
  Let \(I=(f_1,\dots,f_r)\) be a generic ideals in \(\Kxn\), generated
  by forms of degree  
  \(d_1,\dots,d_r\). There is a famous conjecture
  \cite{Moreno:Revlex, Froeberg:Inequality} for the Hilbert
  series of the quotient \(\frac{\Kxn}{I_n}\).
  \begin{conjecture}\label{conj:wk}
    Let \(I=(f_1,\dots,d_r) \subset \Kxn\) be a generic ideal, with
    \(\tdeg{f_i} =d_1\) for \(1 \le i \le r\). Then the Hilbert series
    of the graded algebra \(\frac{\Kxn}{I_n}\) is given by 
    \begin{equation}\label{eq:moreno}
      \left \langle (1-t)^{-n} \prod_{i=1}^r (1-t^{d_i}) \right \rangle
    \end{equation}
  \end{conjecture}
  It is easy to see that if \(r \le n\), the
  generators form a regular sequence, and hence that
  \begin{equation}\label{eq:ci}
    \frac{\Kxn}{I_n}(t) = (1-t)^{-n} \prod_{i=1}^r (1-t^{d_i}) , \qquad
    \text{ for } n \ge r 
  \end{equation}
  In particular, the conjecture holds for  \(r \le n\).
  The conjecture is also know to be true for \(r=n+1\).

  We note that \eqref{eq:ci} implies that 
  \begin{equation}
    \label{eq:analog}
    \lim_{n \to \infty} \frac{ \frac{\Kxn}{I_n}(t)}{\Kxn(t)} =
    \prod_{i=1}^r (1-t^{d_i})
  \end{equation}

  Now suppose that \(I=(f_1,\dots,f_r)\) is a generic ideal in the
  square-free algebra, and that \(f_i\) is a generic form of degree
  \(d_i\). Then 
  \begin{displaymath}
    \frac{\SQ(V_n)}{(f_1,\dots,f_r)} \simeq
    \frac{\Kxn}{(f_1',\dots,f_r',x_1^2,\dots,x_n^2)} 
  \end{displaymath}
  where \(f_i'\) can be taken to be a generic form in \(\Kxn\) which
  maps to \(f_i\) under the canonical epimorphism \(\Kxn
  \twoheadrightarrow \SQ(V_n)\). It seems reasonable to assume that
  the Hilbert series of the quotient will not change if we replace the
  squares of variables with generic quadratic
  forms. Conjecture~\ref{conj:wk} then leads to the following:
  \begin{conjecture}\label{conj:sqfree}
    Let \(r, n, d_1,\dots,d_r\), and let \(I_n\)
    be a generic ideal i \(\SQ(V_n)\) with generators of degrees
    \(d_1,\dots,d_r\). Then 
    \begin{equation}
      \label{eq:sqcase}
      \frac{\SQ(V_n)}{I_n}(t) = \left \langle (1+t)^n \prod_{i=1}^r
        (1-t^{d_i}) \right \rangle
    \end{equation}
  \end{conjecture}
  
  If this conjecture holds (our computations support this), then it
  follows that  
  \begin{equation}
  \label{eq:sqfr}
  \lim_{n \to \infty}\frac{\frac{\SQ(V_n)}{I_n}(t)}{{\SQ(V_n)}(t)} 
  =  
  \prod_{i=1}^r (1-t^{d_i}) 
  \end{equation}
  This is analogous to \eqref{eq:analog}.

\end{section}

\begin{section}{Hilbert series for generic non-principal ideals in the
    exterior algebra} 
  We now throw all caution to the wind to make some bold conjectures
  about the Hilbert series of non-principal generic ideals. 
  Let \(I_n=(f_1,\dots,f_r)\) be a generic ideal in \(\bigwedge V_n\), 
  with \(\tdeg{f_i} = d_i\), and consider the exact sequence
    \begin{equation}
      \label{eq:exseq2}
      0 \longrightarrow \ann(f_r)(-d_r) \longrightarrow
      \frac{\bigwedge V_n}{(f_1,\dots,f_{r-1})}(-d_r)  
      \xrightarrow{\cdot f_r} \frac{\bigwedge
        V_n}{(f_1,\dots,f_{r-1})} \longrightarrow 
      \frac{\bigwedge V_n}{(I)} \longrightarrow 0   
    \end{equation}
    We denote the Hilbert series of \(\frac{\bigwedge V_n}{(I)}\) by
    \(\qHilb\), that of 
    \(\frac{\bigwedge V_n}{(f_1,\dots,f_{r-1})}\) by \(u_n(t)\), and
    that of \(\ann(f_r)\) by  
    \(\aHilb\). Then
    \begin{equation}
      \label{eq:hsum}
      \qHilb = u_n(t) - t^{d_r} u_n(t) + t^d \aHilb.
    \end{equation}
    If \(d_r\) is even, we conjecture
    that \(\aHilb \sim 0\), hence 
    \begin{equation}
      \label{eq:pmeven}
      \qHilb \sim (1-t^{d_r}) u_n(t)
    \end{equation}
    If \(d_r\) is odd, we conjecture that 
    the annihilator of \(f_r\) is ``close'' to the principal ideal on \(f_r\), 
    hence that \(\aHilb \sim (u_n(t) - \qHilb)\), which yields
    \begin{equation}
      \label{eq:pmodd}
      \qHilb (1+t^d) \sim u_n(t)
    \end{equation}
    By induction, we arrive at the following conjecture:
  \begin{equation}\label{eq:hserguessnonp}
    \lim_{n \to \infty} \frac{\qHilb}{(1+t)^n} = \prod_{i=1}^r
    \left(1-(-1)^{d_i} t^{d_i} \right)^{(-1)^{d_i}} \in \Z[[t]],
  \end{equation}
  with respect to the \((t)\)-adic topology.

  One would be tempted to guess that if all \(d_i\)'s are even, the
  Hilbert series of \(\frac{\bigwedge V_n}{(f_1,\dots,f_r)}\) should
  be \emph{exactly}  
  \begin{equation}
    \label{eq:nottrue}
    (1+t)^n \prod_{i=1}^r (1-t^{d_i})
  \end{equation}
  However, this is not true, even for the simplest case \(r=2\) and
  \(d_1=d_2=2\). In Table~\ref{tab:2gens} we tabulate the difference
  between the true Hilbert series and \eqref{eq:nottrue}.

  \tabqgens

\end{section}

\appendix

\begin{section}{The signed incidence matrix has full rank when the
    difference in cardinality is even}

  We prove a ``signed version'' of the well-known theorem that the
  incidence matrix of \(r\)-subsets of \([n]=\set{1,\dots,n}\) into
  \(d+r\)-subsets have full rank. Our proof is a modification of
  the one by Wilson \cite{Wilson:design}.

  To begin, let us define the ``signs'' involved.
  \begin{definition}
    Let \([n]=\set{1,\dots,n}\), and let
    \(\col\) and \(\row\) be two subsets of \([n]\), with
    \begin{displaymath}
      \begin{split}
        \col&=\set{t_1,\dots,t_a}, \qquad t_1 < \dots < t_a \\ 
      \row&=\set{k_1,\dots,k_b}, \qquad k_1 < \dots < k_b 
    \end{split}
  \end{displaymath}
  Then define \(\sigma(\col,\row)\) to be zero if \(\col \not \subseteq \row\), and
  otherwise the 
  sign of the permutation which sorts
  \([\col,\row \setminus \col]\) in ascending order. In other
  words, if \(\col \subseteq \row\) then \(\sigma(\col,\row)\) is the sign of the
  uniquely determined permutation \(\gamma\) such that 
  \begin{displaymath}
    \begin{split}
      t_{\gamma(i)} & = k_i, \qquad 1 \le i \le a \\
      k_{\gamma(j)} & = k_{a+j}, \qquad 1 \le j \le b \\
    \end{split}
  \end{displaymath}
\end{definition}

\begin{definition}\label{def:rowsum}
  Let \([n]=\set{1,\dots,n}\), and let
  \(A\), \(B\) be two subsets of \([n]\), of cardinality \(a\) and
  \(b\), with \(0 \le a < b\). For \(a \le r < b\), we  define
  \begin{equation}
    \label{eq:sdik}
    s_r(A,B,n) = \sum_{\substack{\col \in \binom{[n]}{r}\\
        A \subseteq \col \subseteq B}} \sigma(\col,B)      
    \end{equation}
    For \(0 \le d  \le n\), we define
    \begin{equation}
      \label{eq:sd}
    s_{d,n} = \sum_{R \in \binom{[n]}{d}}  \sigma(R,[n]) 
    = s_d(\emptyset, [n] , n)
    \end{equation}
\end{definition}

\begin{lemma}\label{lemma:sindep}
  With the notations of Definition~\ref{def:rowsum}, put \(d=b-r\).
  We have that
  \begin{equation}
    \label{eq:sindep}
    s_r(A,B,n) = 
    \begin{cases}
      0 & A \not \subseteq B \\
      (-1)^d s_{d,b-a} & A  \subseteq B 
    \end{cases}
  \end{equation}
\end{lemma}
\begin{proof}
  Put \(d=b-r\). If \(A \not \subseteq B\) then clearly
  \(s_r(A,B,n)=0\). Suppose that \(A \subseteq B\).
  Then 
  \begin{displaymath}
    s_r(A,B,n) = \sum_{\substack{\col \in \binom{[n]}{r}\\
        A \subseteq \col \subseteq B}} \sigma(\col,B) 
    = \sum_{\substack{\col \in \binom{B}{r}\\  A \subseteq \col}} 
    \sigma(\col,B),
  \end{displaymath}
  so the sum is independent of \(n\). Furthermore, we can write 
  \(A \subseteq \col \in \binom{B}{r}\) as a disjoint union \(\col = A \cup
  (\col \setminus A)\), hence the sum can be written
  \begin{displaymath}
    \sum_{S \in \binom{B \setminus A}{r-a}} \sigma(S \cup A, B)
     = \sum_{S \in \binom{B \setminus A}{r-a}} \sigma(S , B
      \setminus A).
  \end{displaymath}
  Now, since \(S\) has cardinality \(r-a\), the set \((B \setminus A)
  \setminus S\) has cardinality \(b-a - (r-a) = b-r =d\), so the
  permutation which transforms \([S, B \setminus A]\) to \([B
  \setminus A, S]\) has cardinality \((-1)^d\). Hence, by substituting 
  \(R = (B \setminus A) \setminus S\), we get that the sum is equal to
  \begin{multline*}
    (-1)^d \sum_{S \in \binom{B \setminus A}{v-a}} \sigma((B \setminus
      A) \setminus S , B \setminus A) 
      =  (-1)^d\sum_{R \in \binom{B \setminus A}{d}} \sigma(R, B
      \setminus A) \\ 
      = (-1)^d \sum_{R \in \binom{[b-a]}{d}} \sigma(R, [b-a]),
  \end{multline*}
  which is the desired result.
\end{proof}

\begin{lemma}\label{lemma:notzero}
  Suppose that \(0 < d \le n\), and that \(d\) is even.
  Then \(s_{d,n} > 0\).
\end{lemma}
\begin{proof}
  The lemma is trivially true for \(d=n\). If \(d=2\), we note that
  \(\sigma(\set{v,v+1},[n]) = 1\) for \(1 \le v < n\), since the
  permutation transforming \([v,v+1,1,2,\dots,v-1,v+2,v+3,\dots, n]\)
  to \([1,2,\dots,n]\) is even. Furthermore, the signs of 
  \(\sigma(\set{v,v+\ell},[n])\) alternate in sign as \(\ell\) goes
  from \(1\) to \(n-v\). Thus, for a fixed \(v\), there are either as
  many positive as negative \(\sigma(\set{v,v+\ell},[n])\), or 1 more
  positive than negative, depending on the parity of \(n-v\). By
  summing over all \(v\), we
  conclude that there are always strictly more positive than negative
  signs. 

  Now suppose that we have shown that \(s_{2k',n'} > 0\) for all
  \(k',n'\) such that \(k' < k\). We want to show that that \(s_{2k,n} 
    > 0\). We have that
  \begin{displaymath}
    s_{2k,n} = \sum_{R \in \binom{[n]}{2k}} \sigma(R,[n]),
  \end{displaymath}
  and writing \(R\) as a disjoint union of its first two element, and the
  remaining elements, this becomes
  \begin{multline*}
    \sum_{1 \le k < \ell \le n-2} \sum_{R_2 \in \binom{\set{\ell+1,
          \ell+2,\dots, n}}{2k-2}} \sigma(\set{k,\ell} \cup R_2, [n]) \\
    = \sum_{1 \le k < \ell \le n-2} \sum_{R_2 \in \binom{\set{\ell+1,
          \ell+2,\dots, n}}{2k-2}}\sigma(R_2, \set{\ell+1,
      \ell+2,\dots, n}) 
    = \sum_{1 \le k < \ell \le n-2} s_{2k-2,n-\ell} > 0.
  \end{multline*}
\end{proof}

Next, we define the signed incidence matrix.
\begin{definition}
  Let \(0 < a < b \le n\) be integers. Then \(M_{a,b,n}\) is the
  \(\binom{n}{b}  
  \times \binom{n}{a}\) matrix where the rows are indexed by
  \(b\)-subsets of \([n]\), the columns by \(a\)-subsets of \([n]\),
  and where the entry in row \(B\), column \(A\) is
  \(\sigma(A,B)\). 
\end{definition}

\begin{theorem}\label{thm:fullrank}
  Let \(0 < a < b \le n\) be integers. If \(d=b-a\) is even, then
  \(M_{a,b,n}\) has full rank. 
\end{theorem}
\begin{proof}
  Denote the row indexed by \(\row \in \binom{[n]}{b}\) by \(\tau_\row\), then
  \(\tau_\row\) can be regarded as an element in \(V_a([n])\), the free
  \(\Q\)-vector space on the \(a\)-subsets of \([n]\). If we denote
  the basis element corresponding to a \(a\)-subset \(\col\) by
  \(\epsilon_\col\), then 
  \[\tau_\row = \sum_{\col \in \binom{[n]}{a}}
  \sigma(\col,\row) \epsilon_\col.\]

  The number of rows in \(M_{a,b,n}\) is \(\binom{n}{b}\), and the
  number of columns is  
  \(\binom{n}{a}\). There are less rows than columns if \(a+b > n\),
  as many rows as columns if \(a+b = n\), and more rows than columns
  if \(a+b < n\).

  \begin{enumerate}
  \item If  \(\mathrm{a+b \ge n}\), we must prove that the rows 
  are linearly independent.
  Suppose  that there is a linear relation among the \(\tau_\row\)'s,
  so that 
  \begin{equation}
    \label{eq:linjrel}
    \sum_{\row \in \binom{[n]}{b}} a_\row \tau_\row  = 0
  \end{equation}
  for some numbers \(a_\row\). We shall prove that all \(a_\row = 0\).

  Choose an \(I \subset \binom{[n]}{i}\), \(0 \le i \le a\), and
  define a linear functional \(H_I: V_a([n]) \to \Q\) by
  \begin{equation}
    \label{eq:linfunc}
    f_I(\epsilon_\col) = 
    \begin{cases}
      1 & I \subseteq \col \\
      0 & I \not \subseteq \col 
    \end{cases}
  \end{equation}
  Then if \(\row \in \binom{[n]}{b}\) we have that 
  \begin{multline}
    \label{eq:calc1}
    f_I(\tau_\row) = f_I \left( \sum_{\col \in \binom{[n]}{a}} \sigma(\col,\row)
    \epsilon_\col \right)
    = \sum_{\col \in \binom{[n]}{a}} \sigma(\col,\row) f_I(\epsilon_\col) \\
    = \sum_{I \subseteq \col \subseteq \row} \sigma(\col,\row) 
    = s_a(I,\row,n) = 
    \begin{cases}
      s_{d,b-i} & I \subseteq \row \\
      0 & I \not \subseteq \row \\
    \end{cases}
  \end{multline}
  The last step follows from Lemma~\ref{lemma:sindep}.
  Applying \(f_I\) to \eqref{eq:linjrel} we get that 
  \begin{multline}
    \label{eq:calc2}
      0 = f_I \left( \sum_{\row \in \binom{[n]}{b}} a_\row \tau_\row \right) 
      = \sum_{\row \in \binom{[n]}{b}} a_\row f_I(\tau_\row) \\
      = \sum_{\row \in \binom{[n]}{b}} a_\row s_a(I,\row) 
      = s_{d,b-i} \sum_{\substack{\row \in \binom{[n]}{b} 
          \row \supseteq    I}} a_\row
  \end{multline}
  Since Lemma~\ref{lemma:notzero} tells us that \(s_{d,b-i} \neq 0\), we
  conclude that  
  \begin{equation}
    \label{eq:subsetzero}
    \sum_{\row \supseteq I} a_\row = 0
  \end{equation}

  Now, for any \(J \subset [n]\) we have, by exclusion-inclusion, that 
  \begin{equation}
    \label{eq:inex}
    \sum_{\row \cap J = \emptyset} a_\row = \sum_{I \subset J} (-1)^{\lvert
      I \rvert} \sum_{\row \supseteq I} a_\row
  \end{equation}
  Fix \(\row_0 \in \binom{[n]}{b}\) and put \(J_0 = [n] \setminus
  \row_0\). Since \(\lvert J_0 \rvert = n - b \le a\) we have, using
  \eqref{eq:subsetzero} that  
  \begin{equation}
    \label{eq:a0}
    a_{\row_0} = \sum_{\row \cap J_0 = \emptyset} a_\row = \sum_{I \subseteq
      J_0} (-1)^{\lvert I \rvert} \sum_{\row \supseteq I} a_\row = 0
  \end{equation}
  Since \(a_{\row_0}\) was arbitrary, all \(a_\row\) are zero. This shows
  that the \(\tau_\row\) are linearly independent.

\item If \(\mathrm{n = a+b}\), then \(M\) is a square
  matrix. By the previous case, the vectors \(\tau_\row\) are
  linearly independent, but since there are \(\binom{n}{a} =
  \binom{n}{b}\) such vectors, they form a basis of \(V_a([n])\); in
  particular, they span this vector space.
  
\item  Finally, let us consider the remaining case  \(\mathrm{n > a+b}\), so 
  that there are more rows than columns. We must prove that the rows
  span \(V_a([n])\). We prove this by induction over \(n-a-b\). The
  case \(n-a-b=0\) is already proved, and forms the basis of the
  induction. We assume \(a,b\) fixed, and that the assertion has been
  proved for all \(a+b \le n' < n\).

  Let \(\Gamma \in \binom{[n]}{a}\) be arbitrary. If we can express
  \(\alpha=\epsilon_\Gamma\) as a linear combination of the \(\tau_\row\)'s, we
  are done. To this end, put 
  \begin{equation}
    \label{eq:ap}
    \alpha' = \sum_{\substack{S \in \binom{[n-1]}{a-1}\\ S \cup
        \set{n} = \Gamma}} \epsilon_S \in V_{a-1}([n-1])
  \end{equation}
  Since \(a-1 + b < n-1\), it follows by induction that there are
  scalars \(\setsuchas{d_J}{J \in \binom{[n-1]}{a-1}}\) such that
  \begin{equation}
    \label{eq:Js}
    \alpha' = \sum_{J \in \binom{[n-1]}{a-1}} d_J \tau_J', \qquad
    \tau_J' = \sum_{\substack{S \in \binom{[n-1]}{a-1}\\ S \subseteq
        J}} \epsilon_S
  \end{equation}
  
  For \(\row \in \binom{[n]}{a}\), \(n \in \row\), put \(c_\row' = d_\row \setminus 
  \set{n}\). Define 
  \begin{equation}
    \label{eq:a00}
    \alpha_0 = \sum_{\substack{\row \in \binom{[n]}{a}\\ n \in \row}} 
    c_\row' \tau_\row 
    \in V_a([n]) 
  \end{equation}
  If we write 
  \begin{displaymath}
    \alpha_0 = \sum_{\col \in \binom{[n]}{a}} a_\col' \epsilon_\col
  \end{displaymath}
  we have that for \(\col \in \binom{[n]}{a}\), \(n \in \col\), that 
  \begin{displaymath}
    a_\col' = 
    \begin{cases}
      1 & \col = \Gamma \\
      0 & \col \neq \Gamma
    \end{cases}
  \end{displaymath}
  which implies that 
  \begin{displaymath}
    \alpha_0 = 
    \begin{cases}
      \alpha & n \in \Gamma \\
      0 & n \not \in \Gamma
    \end{cases}
  \end{displaymath}
  In either case, \(\alpha - \alpha_0\) has coordinate 0 in component
  \(\row \in \binom{[n]}{a}\), unless \(n \in \row\). Hence,
  \(\alpha-\alpha_0\) may be regarded as a vector in
  \(V_a([n-1])\). By the induction hypothesis, there exist \(c_\row''\)
  such that
  \begin{equation}
    \label{eq:a0d}
    \alpha-\alpha_0 = \sum_{\row \in \binom{[n-1]}{a}} c_\row'' \tau_\row
  \end{equation}
  Defining 
  \begin{displaymath}
    c_\row = 
    \begin{cases}
      c_\row' & n \in \row \\
      c_\row'' & n \not \in \row
    \end{cases}
  \end{displaymath}
  we get that 
  \begin{displaymath}
    \alpha = \alpha_0 + (\alpha-\alpha_0) =       
    \sum_{\substack{\row \in \binom{[n]}{a}\\ n \in \row}} c_\row' 
    \tau_\row  +
    \sum_{\substack{\row \in \binom{[n]}{a}\\ n \not \in \row}} c_\row'' \tau_\row
    = \sum_{\row \in \binom{[n]}{a}} c_\row \tau_\row
  \end{displaymath}
  \end{enumerate}
\end{proof}

\end{section}

  \begin{section}{Calculations}
    The computer calculations were done on the computers of the UMS
    Medicis, École Polytechnique, and on the computers at the
    Department of Mathematics, Stockholm University. We have used the
    programme Macaulay 2 \cite{MACAULAY2} to calculate
    Hilbert series and 
    minimal free resolutions. To save time and memory, the calculations were
    performed in characteristic 31991. The holes in the tables show
    that there are limits to what we could calculate, even on a
    machine with 2 GB of memory. 
  \end{section}

\bibliographystyle{alpha}
\bibliography{journals,snellman,articles}

\begin{thebibliography}{AHH97}

\bibitem[AAH99]{Aramova:Cohomology}
Anetta Aramova, Luchezar~L. Avramov, and J{\"u}rgen Herzog.
\newblock Resolutions of monomial ideals and cohomology over exterior algebras.
\newblock {\em Transactions of the {A}merican {M}athematical {S}ociety},
  352(2):579--594, 1999.

\bibitem[AHH97]{Aramova:Gotzman}
Annetta Aramova, J{\"u}rgen Herzog, and Takayuki Hibi.
\newblock Gotzman theorems for exterior algebras and combinatorics.
\newblock {\em Journal of {A}lgebra}, 191:174--211, 1997.

\bibitem[Big95]{Bigatti95}
A.~M. Bigatti.
\newblock {\em Aspetti Combinatorici e Computazionali dell'Algebra
  Commutativa}.
\newblock PhD thesis, Universit{\`a} di {T}orino, 1995.

\bibitem[Bou59]{Bourbaki:FormQ}
N.~Bourbaki.
\newblock {\em \'{E}l\'ements de math\'ematique. {P}remi\`ere partie: {L}es
  structures fondamentales de l'analyse. {L}ivre {I}{I}: {A}lg\`ebre.
  {C}hapitre 9: {F}ormes sesquilin\'eaires et formes quadratiques}.
\newblock Hermann, Paris, 1959.
\newblock Actualit\'es Sci. Ind. no. 1272.

\bibitem[CL69]{CleLi}
G.~F. Clements and B.~Lindstr{\"o}m.
\newblock A generalization of a combinatorial theorem of {M}acaulay.
\newblock {\em J. Combinatorial Theory}, 7:230--238, 1969.

\bibitem[Eis95]{Ebud:View}
David Eisenbud.
\newblock {\em Commutative {A}lgebra with a {V}iew {T}oward {A}lgebraic
  {G}eometry}, volume 150 of {\em Graduate {T}exts in {M}athematics}.
\newblock Springer {V}erlag, 1995.

\bibitem[FH94]{Froeberg:HilbGenForm}
Ralf Fr{\"o}berg and Joachim Hollman.
\newblock Hilbert series for ideals generated by generic forms.
\newblock {\em Journal of {S}ymbolic {C}omputation}, 17:149--157, 1994.

\bibitem[FL91]{Froeberg:OnHilb}
Ralf Fr{\"o}berg and Clas L{\"o}fwall.
\newblock On {H}ilbert series for commutative and noncommutative graded
  algebras.
\newblock {\em Journal of {P}ure and {A}pplied {A}lgebra}, 76:33--38, 1991.
\newblock {N}orth {H}olland.

\bibitem[FL00]{KosLie}
Ralf Fr{\"o}berg and Clas L{\"o}fwall.
\newblock Koszul homology and lie algebras with application to generic forms
  and points.
\newblock Technical Report~8, Department of Mathematics, Stockholm University,
  2000.

\bibitem[Fr{\"o}85]{Froeberg:Inequality}
Ralf Fr{\"o}berg.
\newblock An inequality for {H}ilbert series of graded algebras.
\newblock {\em Mathematica {S}candinavica}, 56:117--144, 1985.

\bibitem[GJ73]{Graver:design}
J.~E. Graver and W.~B. Jurkat.
\newblock The module structure of integral designs.
\newblock {\em {J}ournal of {C}ombinatorial {T}heory. {S}eries {A}}, 15:75--90,
  1973.

\bibitem[GS]{MACAULAY2}
Daniel~R. Grayson and Michael~E. Stillman.
\newblock Macaulay 2.
\newblock Computer algebra program, available at
  \texttt{http://www.math.uiuc.edu/Macaulay2/}.

\bibitem[Kan72]{Kantor:incidence}
William~M. Kantor.
\newblock On incidence matrices of finite projective and affine spaces.
\newblock {\em Mathematische {Z}eitschrift}, 124:315--318, 1972.

\bibitem[Kat68]{Katona}
G.~Katona.
\newblock A theorem of finite sets.
\newblock In {\em Theory of graphs (Proc. Colloq., Tihany, 1966)}, pages
  187--207. Academic Press, New York, 1968.

\bibitem[Kru63]{Kruskal}
Joseph~B. Kruskal.
\newblock The number of simplices in a complex.
\newblock In {\em Mathematical optimization techniques}, pages 251--278. Univ.
  of California Press, Berkeley, Calif., 1963.

\bibitem[Mac27]{Macaulay:Enum}
F.~S. Macaulay.
\newblock Some properties of enumeration in the theory of modular systems.
\newblock {\em Proceedings of the London Mathematical Society}, 26:531--555,
  1927.

\bibitem[MS91]{Moreno:Revlex}
Guillermo Moreno-Soc{\'\i}as.
\newblock {\em Autour de la fonction de {H}ilbert-{S}amuel (escaliers
  d'id{\'e}aux polynomiaux)}.
\newblock PhD thesis, {\'E}cole {P}olytechnique, 1991.

\bibitem[Wil73]{Wilson:design}
Richard~M. Wilson.
\newblock The necessary conditions for $t$-designs are sufficient for
  something.
\newblock {\em Utilitas {M}athematica}, 4:207--215, 1973.

\end{thebibliography}

\end{document}